\documentclass[10pt,reqno]{amsart}

\include{psfig}
\usepackage{amssymb,amsmath}
\usepackage{enumerate}
\usepackage{graphicx}
\usepackage{verbatim}
\usepackage{caption}
\usepackage{subcaption}
\usepackage{bcol}
\usepackage[colorinlistoftodos,prependcaption,textsize=small]{todonotes}

\usepackage{natbib}
\newcommand\st{\mbox{s.t. \ }}
\newcommand\model[1]{\hskip -2cm \text{(#1)} \ \ \ \ \ \ \ \ }
\newcommand{\norm}[1]{\ensuremath{\left \|{#1} \right \|}}
\newcommand{\R}{\ensuremath{\mathbb{R}}}

\newtheorem{proposition}{Proposition}

\theoremstyle{definition}

\theoremstyle{remark}

\newcommand{\new}[1]{{\textcolor{black}{#1}}}
\newcommand{\revtwo}[1]{{\textcolor{black}{#1}}}

\title[Conic Integer Programming for Constrained Assortment Optimization]
{A Conic Integer Programming Approach to Constrained Assortment Optimization under the Mixed Multinomial Logit Model}

\author{Alper \c{S}en,  Alper Atamt\"urk and Philip Kaminsky}

\thanks{
	\hskip -6.2mm
	A. \c{S}en: Department
	of Industrial Engineering, Bilkent University Bilkent, Ankara, 06800,
	Turkey \texttt{alpersen@bilkent.edu.tr} \\
	A. Atamt\"urk: Department of Industrial Engineering and Operations
	Research, University of California, Berkeley, CA 94720-1777 USA
	\texttt{atamturk@berkeley.edu} \\
	P. Kaminsky: Department of Industrial Engineering and Operations
	Research, University of California, Berkeley, CA 94720-1777 USA
	\texttt{kaminsky@berkeley.edu}
}

\begin{document}

\maketitle

\begin{abstract}

We consider the constrained assortment optimization problem under the mixed multinomial logit model. Even moderately sized instances of this problem are challenging to solve directly using standard mixed-integer linear optimization formulations. This has motivated recent research exploring customized optimization strategies and approximation techniques.  In contrast, we develop a novel conic quadratic mixed-integer formulation. This new formulation, together with McCormick inequalities exploiting the capacity constraints, enables the solution of large instances using commercial optimization software. \\

\end{abstract}

\begin{center} October 2015; October 2016; May 2017 \end{center}

\BCOLReport{15.06}

\section{Introduction}


\new{Assortment planning, the selection of products that a firm offers to its customers, is a key problem faced by retailers, with direct impact on profitability, market share, and customer satisfaction. A growing stream of operations research literature focuses on assortment optimization problems, where the assortment is optimized to maximize revenue \citep[see][for a review]{kokbookchapter}}.
\ignore{A growing stream of operations research literature focuses on assortment optimization problems, where the selection of products that a firm offers to its customers is optimized to maximize revenue \citep[see][for a review]{kokbookchapter}} To solve this category of problems, customers' purchase behavior must be modeled in a way that captures the impact on the overall demand of product characteristics and customers' substitution between products.  The most commonly used model for customer behavior in this setting is the the multinomial logit (MNL) model, which is based on a probabilistic model of individual customer utilities (see the pioneering work of \citet{ryzin1999relationship}, and followup work by
\citet{cachon2005retail}, \citet{mahajan2001stocking}, \citet{chong2001modeling}, \cite{li2007single}, \cite{rusmevichientong2010dynamic}, \citet{rusmevichientong2012robust}, and
\citet{topaloglu2013joint}). Despite its popularity, the MNL model has several key shortcomings: 1) It relies on the so-called Independence of Irrelevant Alternatives (IIA) assumption, so that a product's market share relative to another product is constant regardless of the other products in the assortment;  2) the total market share of an assortment and the substitution rates within that assortment cannot be independently defined \citep{kok2007demand}. A partial remedy for these problems is possible under an extension of the MNL model, called the nested logit model. Recent work that studies assortment optimization under variants of the nested logit model includes \citet{davis2014assortment}, \citet{gallego2014constrained} and \citet{li2015d}.

In this paper, we consider assortment optimization under a generalization of the MNL model that does not have either of these limitations, the Mixed MNL model (MMNL).  The MMNL model, introduced by \citet{boyd1980effect} and \citet{cardell1980measuring}, has another important characteristic -- as observed by \citet{mcfadden2000mixed}: ``any discrete choice model derived from random utility maximization... can be approximated as closely as one pleases by a MMNL model."   Assortment planning under the MMNL model, also known as the mixtures of MNL model \citep{feldman2015bounding}, MNL with random choice parameters \citep{rusmevichientong2014assortment}, and latent-class MNL \citep{mendez2014branch}, has received considerable interest in the OR/MS community.  The problem also arises as a subproblem in a new approach in revenue management called choice-based deterministic linear optimization that attempts to model customer choice behavior more realistically  \citep{liu2008choice}.

We are particularly interested in assortment optimization under MMNL with constraints on the number of products in the assortment (so-called ``capacity constraints").  While optimal assortments under MNL can be efficiently found \citep{rusmevichientong2010dynamic}, this does not hold true for assortment optimization under MMNL, in either the capacitated or uncapacitated settings.  Indeed, \citet{bront2009} and \citet{rusmevichientong2014assortment} show that the assortment optimization problem under the mixed MNL model is NP-hard.  Motivated by the computational complexity and the ineffectiveness of standard MILP formulations of the problem, \cite{bront2009} propose a greedy heuristic. \citet{mendez2014branch} design and test a branch-and-cut algorithm that generates good but often not provably optimal solutions for both capacitated and uncapacitated versions. \citet{rusmevichientong2014assortment} identify special cases of the (uncapacitated) problem that are polynomially solvable, and characterize the performance of heuristics for other cases. \cite{feldman2015bounding} develop strong upper bounds on the optimal objective value.

In contrast, we show that by formulating this problem in a non-traditional manner, as a conic quadratic mixed-integer program, large instances of the capacitated version of the problem can be solved directly using commercial mathematical optimization software, obviating the need to develop customized heuristics or optimization software to solve the problem.  The advantages of this approach are clear: commercial software is continually developed to take advantages of advances in optimization methods and hardware, it is supported by large software firms, and it allows the inclusion of new constraints without the need for reprogramming. We also show how to further strengthen the formulation with McCormick estimators derived through conditional bounds exploiting the capacity constraints.

\section{Background}

In this section we present a short overview of the mixed multinomial logic model and conic integer optimization.


\subsection{The Consumer Choice Model}
First, recall the traditional MNL model.  Let $N$ be the set of products in the category indexed by $j$. Let $S$ be the assortment -- the subset of products offered by the retailer. Let $\rho_j$ be the unit price for product $j$. The MNL model is based on the utility that a customer gets from consuming a product. For any product, this utility has two components $U_j = u_j + \epsilon_j$, where $u_j$ is a deterministic component and  $\epsilon_j$ is a random component which is assumed to be a Gumbel random variable with mean zero and variance  $\mu^2 \pi^2/6$. Given these, the probability that a customer purchases product $j$ from a given assortment $S$ is  $p_j(S)=\nu_j/(\nu_0+\sum_{k \in S} \nu_k)$,  where $\nu_j=e^{(u_j-\rho_j)/\mu}$ and $\nu_0$ corresponds to the no-purchase option.

As discussed above, we utilize the  mixed MNL (MMNL) model.  This model extends the MNL model by introducing a set $M$ of customer classes. Let $\gamma_i$ be the probability that the demand originates from customer class $i$. The demand in each customer class is governed by a separate MNL model. Let $\nu_{ij}$ be the customer preference for product $j$ in class $i$ and $\nu_{i0}$ be the no-purchase preference in class $i$. Let the unit revenue from product $j$ in class $i$ be $\rho_{ij}$.  We can then write the expected revenue for a given assortment $S$ as
\begin{equation}
\sum_{i \in M} \gamma_i \left[ \frac{\sum_{j \in S}  \rho_{ij} \nu_{ij}}{\nu_{i0}+\sum_{j \in S} \nu_{ij}} \right].
\end{equation}
Due to space or administrative restrictions, there can be various constraints on the depth of the assortment that can be carried. Let $K$ be the set of resources that may constrain the assortment. Let $\beta_{kj}$ denote the amount of resource $k$ used by product $j$ and let $\kappa_k$ denote the amount of resource $k$ available.
The capacitated assortment optimization problem is therefore to select the assortment $S$ in this setting. 

\subsection{Conic Integer Optimization}

\label{sec:conicmip}

Conic optimization refers to optimization of a linear function over
conic inequalities \citep{N:convex-book}. 
A conic quadratic constraint on $x \in \mathbf{R}^n$ is a constraint of the form
\[
 \norm{A x - b} \ \le c' x - d .
\]
Here \revtwo{$\norm{\cdot}$} is the L2 norm, $A$ is an $m\times n$-matrix, $b$ is an $m$-column vector, $c$ is an $n$-column vector,
and $d$ is a scalar. We refer the reader to \cite{LBL:socp-apps} and \cite{AG:socp} for reviews of conic quadratic optimization and its applications.

Although there is an extensive body of literature on convex conic quadratic optimization,
development of conic optimization with integer variables is quite recent 
\citep{CI:conicmip, AN:conicmir-ipco, AN:coniclift,AMP:conicknap-gub}.
With the growing availability of commercial solvers for these problems (for example, both CPLEX and Gurobi now include solvers for these models), conic quadratic integer models have recently been employed to address problems in portfolio optimization \citep{VAN:cmip-bb}, value-at-risk minimization \citep{AN:conicobj}, machine scheduling \citep{AAG:matchup}, and supply chain network design \citep{ABS:riskpool}, airline rescheduling with speed control\citep{AAG:cruise}.
However, to the best of our knowledge, this approach has not been previously used to solve assortment optimization problems.

Conic quadratic inequalities are often used to represent
a rotated cone/hyperbolic inequality
\begin{equation} \label{eq:quad-cone1}
x_1^2 \le x_2 x_3,
\end{equation}
on $x_1, x_2, x_3 \ge 0$. It is easily verified that
hyperbolic inequality~\eqref{eq:quad-cone1}
can then be equivalently written as a conic quadratic inequality
\begin{equation} \label{c:cq1}
\norm{(2x_1, x_2-x_3)} \le x_2+x_3.
\end{equation}
In our conic reformulation of the assortment optimization problem, we make use of \revtwo{the rotated cone} inequalities \revtwo{\eqref{eq:quad-cone1} in our models}.

\section{The Capacitated Assortment Optimization Problem}

In this section we first recall the traditional MILP formulation of the capacitated assortment optimization problem, and then present an alternative conic quadratic mixed 0-1 formulation of the problem, and strengthen the formulation using McCormick estimators based on conditional bounds.

\subsection{The Traditional MILP Formulation}
Given the MMNL demand model, define $x_j$ to be 1 if product $j$ is offered in the assortment and 0 otherwise.  We can then state the capacitated assortment optimization problem as a nonlinear binary optimization:
\begin{align}
\label{eq:modelfirst}
\max         &  \sum_{i \in M} \gamma_i \left[ \frac{\sum_{j \in N}  \rho_{ij} \nu_{ij} x_j}{\nu_{i0}+\sum_{j \in N} \nu_{ij}x_j} \right] \\
\model{CAOP} \st  & \new{\sum_{j \in N} \beta_{kj} x_{j}  \leq   \kappa_{k},\; \; \forall k \in K } \\
\label{eq:modellast}  &    x_j \in \{0,1\},\; \; \forall j \in N.
\end{align}

Traditionally, (CAOP) is formulated as a Mixed Integer Linear Program \citep[see, e.g.][]{bront2009,mendez2014branch}.
First, letting $y_i=1/(\nu_{i0}+\sum_{j \in N} \nu_{ij}x_j)$, the problem can be posed as a bilinear mixed 0-1 optimization problem:
\begin{align}
\label{eq:modelfirstMILP1}
\max  &  \sum_{i \in M}  \sum_{j \in N} \gamma_i  \rho_{ij} \nu_{ij} y_i x_j \\
\st &  \new{ \sum_{j \in N} \beta_{kj} x_{j}  \leq   \kappa_{k},  \; \; \forall k \in K }\\
\model{CAOP'} & \nu_{i0} y_i + \sum_{j \in N} \nu_{ij} y_i x_j  =  1, \; \; \forall i \in M \\
                            &  y_i  \geq  0, \; \; \forall i \in M \\
\label{eq:modellastMILP1}   &  x_j \in \{0,1\},\; \; \forall j \in N.
\end{align}
The bilinear terms $y_i x_j$ in the formulation can be linearized using the standard ``big-$M$" approach: For any bilinear term $y x$, where $y$ is continuous and non-negative and $x$ is binary, define a new continuous variable $z=yx$ and add the following inequalities to the formulation: $y-z\leq U(1-x)$, $0 \le z\leq y$ and $z \leq U x$, where $U$ is a sufficiently {\it large upper bound on $y$}.  Employing this technique, and selecting $1/\nu_{i0}$ for $U$, leads to the following mixed-integer linear formulation:
\begin{align}
\label{eq:modelfirstMILP2}
\max        &  \sum_{i \in M} \sum_{j \in N} \gamma_i  \rho_{ij} \nu_{ij} z_{ij} \\
\label{eq:capacity}
\st &  \new{\sum_{j \in N} \beta_{kj} x_{j}  \leq   \kappa_{k}, \; \; \forall k \in K } \\
& \nu_{i0} y_i + \sum_{j \in N} \nu_{ij} z_{ij}  =  1, \; \; \forall i \in M \\
\model{MILP}
\label{MILP:globalbounds1}
& \nu_{i0} (y_i-z_{ij})   \leq  1-x_j, \;\; \forall i \in M, \; \; \forall j \in N \\
& 0 \le z_{ij}  \leq  y_i,  \;\; \forall i \in M \; \; \forall j \in N,\\
\label{MILP:globalbounds3}
& \nu_{i0} z_{ij}  \leq  x_j,    \;\; \forall i \in M, \; \; \forall j \in N \\
  &   x_j \in \{0,1\},\; \; \forall j \in N  \\
& z_{ij}  \geq  0,   \; \; i \in M \; \; j \in N \\
\label{eq:modellastMILP2} & y_i  \geq  0, \; \; i \in M.
\end{align}

As shown in \citet{bront2009,mendez2014branch,feldman2015bounding}, formulation (MILP) does not scale well.
In particular, when the capacity constraints (\ref{eq:capacity}) are tight, solution times are prohibitive even
for moderately sized instances.

\subsection{The Conic Formulation}

In order to give a conic reformulation, we first restate the objective as minimization.
Letting $\overline{\rho}_i=\max_{j \in N} \rho_{ij}$,
the objective (\ref{eq:modelfirst}) of (CAOP) can be written as
\begin{equation}
\label{obj_rewritten}
\max \sum_{i \in M} \gamma_i \overline{\rho}_i  - \sum_{i \in M} \gamma_i \left[ \frac{\nu_{i0} \overline{\rho}_i+\sum_{j \in N} \nu_{ij} (\overline{\rho}_i-\rho_{ij}) x_j}{\nu_{i0}+ \sum_{j \in N} \nu_{ij} x_j} \right].
\end{equation}
As the first component  in (\ref{obj_rewritten}) is constant, we can pose the problem as minimizing the second component in (\ref{obj_rewritten}). Also, since the objective coefficients are nonnegative, it suffices to use only lower bounds on $y$ and $z$ variables, leading to:
\begin{align}
\label{min:objective1}
\min        &   \sum_{i \in M} \gamma_i \nu_{i0} \overline{\rho}_i y_i +  \sum_{i \in M} \sum_{j \in N} \gamma_i \nu_{ij} (\overline{\rho}_i-\rho_{ij}) z_{ij} \\
\st & \new{ \sum_{j \in N} \beta_{kj} x_{j}  \leq   \kappa_{k}, \;\; k \in K } \\
\label{eq:prereform}  \model{CAOP''}
       &                   z_{ij} \ge y_i x_j, \;  \; i \in M, \; \; j \in N  \\
\label{eq:prereform2}
       &                   y_i \geq \frac{1}{\nu_{i0}+ \sum_{j \in N} \nu_{ij} x_j},  \; \; i \in M  \\
    & x_j \in \{0,1\},\; \; \; j \in N \\
    & z_{ij}   \geq 0, \; \; \; \; i \in M, \; \; j \in N \\
    \label{eq:facet2}
    & y_i 		\geq 0, \; \; \; \; i \in M.				
\end{align}

\ignore{
In order to represent (\ref{eq:prereform}) and  (\ref{eq:prereform2}) as conic constraints, first define
\begin{equation}
w_i={\nu_{i0}+ \sum_{j \in N} \nu_{ij} x_j}.
\end{equation}
Then the constraints (\ref{eq:prereform}) and (\ref{eq:prereform2}) can be written as
\begin{eqnarray}
w_i+z_{ij} &\geq& \sqrt{(w_i-z_{ij})^2+(2x_j)^2},\\
w_i+y_i &\geq& \sqrt{(w_i-y_i)^2+2^2}.
\end{eqnarray}

The resulting inequalities are known as rotated cone constraints:
}

\ignore{
\noindent
Now defining
\begin{equation}
w_i={\nu_{i0}+ \sum_{j \in N} \nu_{ij} x_j}
\end{equation}
and {\new{red} observing that $w_i \geq 0$ for all $i \in M$} and $x_j = x_j^2$, we can state constraints (\ref{eq:prereform}) and  (\ref{eq:prereform2})
in rotated cone form:
\begin{eqnarray}
z_{ij} w_i &\geq& { x_j^2} \\
y_i w_i &\geq& 1.
\end{eqnarray}
}

\noindent
\new{Observe that constraints \eqref{eq:prereform}--\eqref{eq:prereform2} are satisfied at equality at an optimal solution.}
Now, defining
\begin{equation}
w_i={\nu_{i0}+ \sum_{j \in N} \nu_{ij} x_j}
\end{equation}
and observing that $w \geq 0$,
one can state constraints (\ref{eq:prereform2}) in rotated cone form:
\begin{eqnarray}
\label{eq:ywis1}
y_i w_i &\geq& 1.
\end{eqnarray}

\noindent
\new{As constraint \eqref{eq:ywis1} is satisfied at equality at an optimal solution, $w \ge 0$,} and $x_j = x_j^2$ \new{for a binary vector $x$}, constraint (\ref{eq:prereform}) can also be stated in rotated cone form:
\begin{eqnarray}
z_{ij} w_i &\geq& { x_j^2}
\end{eqnarray}

\noindent
Although redundant for the mixed-integer formulation, we also use the constraints
\begin{align} \label{eq:facet}
& \nu_{i0} y_i + \sum_{j \in N} \nu_{ij} z_{ij}  \ge  1, \; \; \forall i \in M
\end{align}
to strengthen the continuous relaxation of the formulation.

The final conic quadratic mixed 0-1 program is, therefore:
\begin{align}
\label{eq:finalconicfirst}
\min        &   \sum_{i \in M} \gamma_i \nu_{i0} \overline{\rho}_i y_i +  \sum_{i \in M} \sum_{j \in N} \gamma_i \nu_{ij} (\overline{\rho}_i-\rho_{ij}) z_{ij} \\
\st & \new{  \sum_{j \in N} \beta_{kj} x_{j}  \leq   \kappa_{k}, \; \; k \in K } \\
\label{eq:w}
& w_i  =  {\nu_{i0}+ \sum_{j \in N} \nu_{ij} x_j} ,  \; \; i \in M \\
\label{eq:tij}
\model{CONIC} & z_{ij} w_i\geq  x_j^2,  \; \; i \in M \; \; j \in N \\
\label{eq:yw}
& y_i w_i \geq 1,  \; \; i \in M \\
& \nu_{i0} y_i + \sum_{j \in N} \nu_{ij} z_{ij}  \ge  1, \; \; \forall i \in M \\
& x_j \in \{0,1\},\; \; j \in N \\
& z_{ij}  \geq  0,   \; \; i \in M, \; \; j \in N \\
\label{eq:finalconiclast}
& y_i  \geq  0, \; \; i \in M.
\end{align}

In contrast to the traditional formulation (MILP), the conic formulation does not require ``big-$M$" constants for linearization,
which lead to weak LP relaxations especially for the tightly capacitated cases. On the other hand, for the conic formulation when capacity is low, small values of $w_i$ tighten the constraints $z_{ij} w_i\geq  x_j^2$, leading to stronger bounds. \revtwo{The next proposition provides a theoretical justification for adding inequalities \eqref{eq:facet}
to the formulation. Preliminary computations also showed a significant strengthening of the conic formulation with the addition of inequalities \eqref{eq:facet}.}

\begin{proposition}
		\label{prop:facet} Inequality \eqref{eq:facet} is a facet-defining for $\text{\revtwo{cl} conv}\{(x,y,z) \in \{0,1\}^N \times \R^M \times \R^{M\times N}: \eqref{eq:prereform}-\eqref{eq:facet2} \}.$
	\end{proposition}
	Proof.
	Let $S = \{(x,y,z) \in \{0,1\}^N \times \R^M \times \R^{M\times N}: \eqref{eq:prereform}-\eqref{eq:facet2} \}.$
	First, observe that even though constraints \eqref{eq:prereform2} are nonlinear, $S$ is a union of polyhedra (one polyhedron for each assignment of the binary variables); hence, $\text{\revtwo{cl} conv}(S)$ is a polyhedron. Let $e_k$ be the $k$th unit vector, $\hat y = \sum_{k \in M}  e_k/\nu_{k0}$ and
	$\revtwo{\hat y^i} = \sum_{k \in M \setminus\{i\}}  e_k/\nu_{k0}$.
	Consider the following $|N|+|M|+|M||N|$  points in $S$ satisfying $\nu_{i0} y_i + \sum_{j \in N} \nu_{ij} z_{ij}  =  1$: $ (0,  \hat y , 0);
	(0,  \hat y + \epsilon e_k, 0), \ k \in M \setminus \{i\}, \ \epsilon > 0;
	(0, \hat y, \epsilon e_{kj}), \ k \in M \setminus \{i\}, j \in N, \ \revtwo{0 < \epsilon < 1};
	(e_j, \revtwo{\hat y^i} + e_i/(\nu_{i0}+\nu_{ij}), e_{ij}/(\nu_{i0} + \nu_{ij}), \ j \in N;
	(e_j, \revtwo{\hat y^i} + (1 - \epsilon)e_i/(\nu_{i0} + \nu_{ij}), (1 + (\nu_{i0} \epsilon/\nu_{ij}))e_{ij}/(\nu_{i0} + \nu_{ij})), j \in N, \ 0 < \epsilon < 1$, 
	\revtwo{where $e_{ij}$ is the $ij$th unit vector}.
	It is easily checked that these points are affinely independent.

\subsection{McCormick Estimators}
\label{sec:mc}

The capacitated assortment formulations can be further strengthened using McCormick
estimators for the bilinear terms. To that end,
we give simple upper and lower bounds on
\begin{equation}
y_i= \frac{1}{\nu_{i0}+\sum_{j \in N} \nu_{ij}x_j}, \ i \in M.
\end{equation}
\new{
The lower bounds make use of the capacity constraints \eqref{eq:capacity}. For $i \in M$,  define the auxiliary problem
\begin{align}
\label{multipleknapsack}
f_i=\max        &   \sum_{j \in N} \nu_{ij} x_j \\
\model{BND} \ \st &   \sum_{j \in N} \beta_{kj} x_{j}  \leq   \kappa_{k}, \; \; k \in K \\
\label{knapsacklast}
& x_j \in \{0,1\},\; \; j \in N.
\end{align}
\begin{proposition} \label{prop:bounds} The following bounds on variables $y_i$, $i \in M$, are valid:
\begin{eqnarray}
y_i^{\ell} &:=& \frac{1}{\nu_{i0}+f_i} \le y_i \\
y_i^{u}&:=&\frac{1}{\nu_{i0}} \ge y_i.
\end{eqnarray}
\end{proposition}
Proposition~\ref{prop:bounds} provides global bounds on variables $y$. Next, we give conditional bounds.
Let $f_{i|x_j=\xi}$ be the objective function value of (BND) when an additional constraint $x_j=\xi, \ j \in N$ is imposed. 
\begin{proposition} \label{prop:bounds-cond} For $j \in N$, the following conditional bounds on variables $y_i$, $i \in M$, are valid:
\begin{align}
x_j = 0 \ \Rightarrow & \
\begin{aligned}
\ \ \ y_{i|x_j=0}^\ell := \frac{1}{\nu_{i0} + f_{i|x_j=0}}  \le y_i
\end{aligned} \\
x_j = 1 \ \Rightarrow & \ \left \{
\begin{aligned}
 y_{i|x_j=1}^\ell := & \frac{1}{\nu_{i0} + f_{i|x_j=1}}  \le y_i \\
 y_{i|x_j=1}^u := & \frac{1}{\nu_{i0} + \nu_{ij}} \ge y_i.
 \end{aligned}
 \right .
\end{align}
\end{proposition}
(BND) is a binary multiple constraint knapsack problem, so it may be prohibitive to find the optimal $f_i$ and $f_{i|x_j=\xi}$  except in special cases. However, note that to get a lower bound on $y_i$, an upper bound on the optimal value of (BND) is sufficient, and this can be found by solving an easier relaxation of the problem, e.g., the linear optimization relaxation.
\ignore{
For a looser bound, one can consider each constraint separately. One can then use a simple greedy algorithm that fills the knapsack starting with products with higher $\nu_{ij}/\beta_{kj}$ ratios for each problem and pick the best bound.}
}

For the special case of a single cardinality constraint, one can obtain exact closed form lower bounds on $y$.

\begin{proposition} \label{prop:boundssinglecardinality} For a single cardinality constraint of the form $\sum_{j \in N} x_j \leq \kappa$, the following global and conditional lower bounds on $y_i, i \in M$, are valid:
\begin{eqnarray}
y_i^{\ell} &:=& \frac{1}{\nu_{i0}+\sum_{k=1}^{\kappa} {\nu_{i{[k]}}}}  \\
y_{i|x_j=0}^\ell &:=& \frac{1}{\nu_{i0} + \sum_{k=1}^{\kappa} {\bar \nu_{i{[k]}}}} \\
y_{i|x_j=1}^\ell &:= & \frac{1}{\nu_{i0} + \nu_{ij} + \sum_{k=1}^{\kappa-1} {\bar \nu_{i{[k]}}}},
\end{eqnarray}
where $\nu_{i{[k]}}$ is defined as the $k$th largest of preferences $\nu_{im}, \, m \in N$ and $\bar \nu_{i{[k]}}$ is defined as the $k$th largest of preferences $\nu_{im}, \, m \in N \setminus \{j\}$.
\end{proposition}
\new{Similar exact closed-form bounds can be developed
when there are multiple non-overlapping cardinality constraints (i.e., the assortment can contain at most a fixed number of products from each product sub-group).}

Using the global and conditional bounds on $y_i$, $i \in M$, above, one can write the following valid  McCormick inequalities \citep{mc:factorable}
for each bilinear term $z_{ij} = y_i x_j$: 
\begin{eqnarray}
\label{eq:also_used_elsewhere}
			        z_{ij} &\leq& y_{i|x_j=1}^u x_j, \; \;  i \in M, \; j \in N \label{mc1s} \\		
\model{MC}		  z_{ij} &\geq& y_{i|x_j=1}^\ell x_j, \; i \in M, \; j \in N  \label{mc4s} \\
						  z_{ij} &\leq& y_i - y_{i|x_j=0}^{\ell} (1-x_j), \;  i \in M, \; j \in N   \label{mc2s} \\
						  z_{ij} &\geq& y_i - y_i^{u} (1-x_j),  \; i \in M, \; j \in N.  \label{mc3} 						
\end{eqnarray}

Note that the inequality (\ref{eq:also_used_elsewhere}) is also used in \citet{mendez2014branch} and that (\ref{mc3}) is the same as (\ref{MILP:globalbounds1}) in model (MILP).

\new{
Based on the discussion thus far,
four different formulations can be used to solve the capacitated assortment optimization problem under MMNL. The first one is (MILP), which
can be strengthened by replacing constraints \eqref{MILP:globalbounds1}-\eqref{MILP:globalbounds3} with the stronger McCormick estimators (MC). We denote this strengthened formulation  (MILP+MC).
The third formulation is (CONIC), which can also be strengthened by adding McCormick inequalities (MC). This fourth formulation is denoted (CONIC+MC). Notice that one can convert (MILP) and (MILP+MC) to minimization problems by using the equivalent objective (\ref{eq:finalconicfirst}). This leads to the observation that (CONIC+MC) is a strengthening of (MILP+MC) with constraints (\ref{eq:w}), (\ref{eq:tij}), and (\ref{eq:yw}). Therefore, (CONIC+MC) is stronger than (MILP+MC), which is itself stronger than (MILP). The numerical experiments reported in the next section show the significance and the effect of differences in the strength of these formulations.
}

\section{Numerical Study}

In order to test the effectiveness of the conic optimization approach and the McCormick inequalities, we perform a numerical study on four sets of problems. \new{The optimization problems are solved with Gurobi 6.5.1 solver on a computer with an Intel Core i7-4510U 2.00 GHz (2.60 GHz Turbo) processor and 8 GB RAM operating on 64-bit Windows 10. We use the default settings of Gurobi except that we force the solver to use the linear outer-approximation approach when solving continuous relaxations of conic programs. The outer-approximation allows warm starts with the dual simplex method and speeds up solving node relaxations.} The time limit is set to 600 seconds.

The first set of problems are created by randomly generating instances with $|N|=200$ products and $|M|=20$ customer classes . The product prices are the same across the customer classes ($\rho_{ij}=\rho_j$) and are drawn from a uniform[1,3] distribution. The preferences $\rho_{ij}$ are drawn from a uniform[0,1] distribution.  The parameter $\gamma_i=1/20$ for all $i \in M$.  The  no purchase parameter $\nu_{i0}=\nu_{0}$ is either 5 or 10. The capacity constraint is in the form of a cardinality constraint. The maximum cardinality $\kappa$ of the assortment is one of five possible values: $\{10,20,50,100,200\}$.  For each of these $5 \times 2 = 10$ capacity and no-purchase probability combinations, we generate five instances, resulting in a total of 50 instances. All data files are available at \texttt{http://ieor.berkeley.edu/$\sim$atamturk/data/assortment.optimization} .

\new{
We test the effectiveness of four formulations: (MILP), (MILP+MC), (CONIC) and (CONIC+MC). In addition, we compare these with the formulation of \cite{mendez2014branch}, which strengthen (MILP) by replacing (\ref{MILP:globalbounds3}) with (\ref{eq:also_used_elsewhere}) and
by introducing five classes of valid inequalities.
Three of these are polynomial in the size of the model, while the rest are exponential.
We run their formulation using the three classes of polynomial valid inequalities.
}

\bgroup
\setlength\tabcolsep{0.08cm}
\def\arraystretch{1}%
\begin{table}[h]
\caption{Results for problems with 200 products and 20 classes.}
\label{smaller_final}
\begin{center}
\begin{footnotesize}
\begin{tabular}{c|c|r|rc|rc|rc|rc|rc} \hline \hline
\rule{0mm}{1.2em} &  &   & \multicolumn{2}{c|}{MILP} & \multicolumn{2}{c|}{MILP+MC} & \multicolumn{2}{c|}{M\'{e}ndez-D\'{i}az} & \multicolumn{2}{c|}{CONIC} & \multicolumn{2}{c}{CONIC+MC} \\
\cline{4-13}
$\nu_0$ & $\kappa$ & \revtwo{assort} &\multicolumn{1}{c}{rgap} & \multicolumn{1}{c|}{time/\#} &  \multicolumn{1}{c}{rgap} & \multicolumn{1}{c|}{time/\#} &\multicolumn{1}{c}{rgap} & \multicolumn{1}{c|}{time/\#} &\multicolumn{1}{c}{rgap} & \multicolumn{1}{c|}{time/\#} &\multicolumn{1}{c}{rgap} & \multicolumn{1}{c}{time/\#} \\
& & \revtwo{bind} & \multicolumn{1}{c}{egap} & \multicolumn{1}{c|}{nodes} & \multicolumn{1}{c}{egap} & \multicolumn{1}{c|}{nodes}& \multicolumn{1}{c}{egap} & \multicolumn{1}{c|}{nodes}& \multicolumn{1}{c}{egap} & \multicolumn{1}{c|}{nodes}& \multicolumn{1}{c}{egap} & \multicolumn{1}{c}{nodes} \\
\hline \hline
 	&	10	& \revtwo{10.0} &	52.56\%	&	--	&	12.33\%	&	--	&	51.46\%	&	--	&	3.20\%	&	32.82/5	&	0.27\%	&	8.72/5	\\
	&		  & \revtwo{5}  &	45.10\%	&	3076	&	9.85\%	&	6374	&	50.16\%	&	0	&	0.00\%	&	1449	&	0.00\%	&	14	\\
\cline{2-13}
	&	20	& \revtwo{20.0} &	33.38\%	&	--	&	10.25\%	&	--	&	33.37\%	&	--	&	5.88\%	&	122.74/4	&	0.36\%	&	9.58/5	\\
	&		  & \revtwo{5}  &	32.07\%	&	11626	&	8.34\%	&	13819	&	33.36\%	&	44.8	&	0.10\%	&	2851	&	0.00\%	&	23	\\
\cline{2-13}
5 	  &	50 & \revtwo{50.0}	&	2.81\%	&	481.16/2	&	0.94\%	&	27.73/3	&	2.79\%	&	--	&	17.14\%	&	--	&	0.02\%	&	2.38/5	\\
    	&		 & \revtwo{5}  &	1.72\%	&	87695	&	0.09\%	&	27779	&	2.78\%	&	102.6	&	2.94\%	&	1566	&	0.00\%	&	0	\\
\cline{2-13}
 	&	100	& \revtwo{65.4}&	0.08\%	&	4.26/5	&	0.03\%	&	1.22/5	&	0.07\%	&	366.16/5	&	23.66\%	&	--	&	0.01\%	&	1.82/5	\\
 	&		  & \revtwo{0} &	0.00\%	&	790	&	0.00\%	&	0	&	0.00\%	&	124	&	7.23\%	&	768	&	0.00\%	&	0	\\
\cline{2-13}
	&	200	& \revtwo{65.4}&	0.08\%	&	2.29/5	&	0.04\%	&	1.06/5	&	0.07\%	&	366.57/5	&	23.66\%	&	--	&	0.01\%	&	1.92/5	\\
	&		  & \revtwo{0} &	0.00\%	&	343	&	0.00\%	&	0	&	0.00\%	&	117.6	&	13.12\%	&	747	&	0.00\%	&	0	\\
\hline	
	&	10	& \revtwo{10.0} &	24.74\%	&	--	&	7.20\%	&	--	&	20.69\%	&	--	&	1.93\%	&	22.50/5	&	0.10\%	&	6.47/5	\\
	&		  & \revtwo{5} &	10.26\%	&	47690	&	5.44\%	&	6555	&	19.70\%	&	0.2	&	0.00\%	&	1054	&	0.00\%	&	4	\\
\cline{2-13}
	&	20	& \revtwo{20.0} &	38.66\%	&	--	&	8.47\%	&	--	&	38.65\%	&	--	&	3.61\%	&	86.77/5	&	0.16\%	&	8.62/5	\\
	&		  & \revtwo{5} &	31.57\%	&	1613	&	7.20\%	&	9498	&	38.61\%	&	4.2	&	0.00\%	&	1374	&	0.00\%	&	7	\\
\cline{2-13}
10 		&	50 & \revtwo{50.0} &	10.50\%	&	--	&	2.92\%	&	--	&	10.50\%	&	--	&	10.31\%	&	--	&	0.08\%	&	7.37/5	\\
			&		 & \revtwo{5} & 9.89\%	&	25276	&	2.02\%	&	31281	&	10.49\%	&	48.6	&	1.30\%	&	1454	&	0.00\%	&	72	\\
\cline{2-13}
 	&	100	& \revtwo{91.8}	& 0.04\%	&	3.46/5	&	0.01\%	&	1.20/5	&	0.03\%	&	306.05/5	&	18.40\%	&	--	&	0.00\%	&	1.77/5	\\
 	&			&	\revtwo{1} & 0.00\%	&	406	&	0.00\%	&	0	&	0.00\%	&	255.8	&	4.62\%	&	766	&	0.00\%	&	0	\\
\cline{2-13}
 	&	200	&	\revtwo{92.0} & 0.04\%	&	2.89/5	&	0.01\%	&	0.93/5	&	0.03\%	&	282.31/5	&	18.41\%	&	--	&	0.00\%	&	1.67/5	\\
	&			&	\revtwo{0} & 0.00\%	&	462	&	0.00\%	&	0	&	0.00\%	&	82.4	&	5.86\%	&	768	&	0.00\%	&	0	\\
\hline
\multicolumn{3}{c|}{Average} 	 &	16.29\%	&	46.67/22	&	4.22\%	&	4.58/23	&	15.76\%	&	330.27/20	&	12.62\%	&	63.23/19	&	0.10\%	&	5.03/50	\\
\multicolumn{3}{c|}{}				& 13.06\%	&	17898	&	3.30\%	&	9531	&	15.51\%	&	78.02	&	3.52\%	&	1280	&	0.00\%	&	12	\\
\hline \hline	
\end{tabular}
\end{footnotesize}
\end{center}
\end{table}
\egroup

\new{
Table \ref{smaller_final} presents averages of root gap, end gap, solution time and the number of search nodes over five instances for each no purchase preference $\nu_0$, capacity level $\kappa$ and formulation. \revtwo{The number of products in the assortment ($\sum_{j \in N} x_j^*$, averaged over five instances) and the number of instances where the capacity is binding in the optimal solution are given by assort and bind, respectively.} The root gap is computed as $\text{rgap} = 100 \times \revtwo{(\text{zopt} - \text{zroot}) |/|\text{zopt}|}$, where zroot is the objective value of the continuous relaxation \revtwo{(before presolve and root cuts)} and zopt is the value of the optimal integer solution.
The end gap is computed as $\text{egap} = 100 \times \revtwo{(\text{zopt} - \text{zbb})|/|\text{zopt}|}$, where zbb is the best lower bound at termination. If an instance is solved to optimality zbb equals zopt \revtwo{(within the default optimality gap 0.01)}.
In the tables, time refers to the average solution time (in seconds) for the instances that are solved within the time limit and \# refers to
the number of instances solved within the time limit.
The last row reports the averages for rgap, egap, time and nodes and the total number of instances solved.
}


\new{
As observed in previous studies, the traditional (MILP) formulation performs poorly, except when the capacity constraint is loose. The time limit is reached for 28 instances with tight capacity constraints. The poor performance appears to be due to the weak relaxation, leading to excessive branching. The remaining gaps at termination are quite large for the unsolved instances.
With the addition of McCormick inequalities (MC), root and end gaps improve substantially in all cases. The average root gap drops from 16.29\% to 4.22\%. However, this is still not enough to solve the capacitated cases. McCormick inequalities help to solve only one additional instance within the time limit.
}

\new{
For our data set, the polynomial inequalities of \cite{mendez2014branch} lead to \revtwo{a small} reduction in root gaps compared to (MILP). 
\revtwo{Cutting plane algorithms implementing separation for the exponential classes of inequalities of \cite{mendez2014branch} may lead to a further reduction. Although we use a different data set, consistent with their numerical study, the M\'endez-D\'iaz formulation is more effective for high capacity instances.} 
Model (MILP+MC) is considerably stronger. The strength of (MILP+MC) over M\'{e}ndez-D\'{i}az is due to conditional McCormick inequalities \eqref{mc4s} and \eqref{mc2s} based on strong lower bounds on $y$.
}


In contrast to the linear formulations, most of the capacitated instances are solved easily with the conic formulation. This is due to small root gaps, leading to only limited enumeration. However, the performance of the conic formulation degrades for high capacity instances. \revtwo{Observe that MILP and CONIC formulations are not directly comparable. The CONIC formulation may be weaker than the MILP formulation for high capacity instances; whereas the MILP formulation tends to be weaker than the CONIC formulation for low capacity instances.}

\new{
The results are dramatically better when the McCormick inequalities are added to the conic formulation. The average root gap drops to a mere 0.10\% and \textit{all} instances are solved to optimality, on average, in five seconds. On average, only 12 nodes are needed in the search tree. For some instances, the CONIC+MC is more than 100 times faster than the other approaches. This is due to the joint effect of the tightening of the formulation using conic constraints and McCormick inequalities as observed with very small root gaps for all instances. \revtwo{As noted in Section~\ref{sec:mc} CONIC+MC dominates MILP+MC.}
}

In Table \ref{larger_final}, we report the results of experiments for instances with 500 products and 50 classes. The preference values and prices are generated as before.  Each class again has equal weight ($\gamma_i=1/50$). The capacity $\kappa$ is one of $\{20,50,100,200,500\}$, and the no purchase parameter $\nu_{i0}$ is either 10 or 20. \new{Since \revtwo{our experiments do not indicate a significant improvement from employing} the approach in \cite{mendez2014branch} over (MILP), \ignore{and performs worse than (MILP+MC) for the smaller instances} we do not include it for the remaining experiments. We also note that five instances cannot be solved using any of the formulations within the time limit. For those instances, the optimal integer solutions are obtained separately using CONIC+MC formulation by extending the time limit. Therefore, root gap and end gap calculations are still with respect to the optimal integer solutions.}
\bgroup
\setlength\tabcolsep{0.08cm}
\def\arraystretch{1}%
\begin{table}[h]
\caption{Results for problems with 500 products and 50 classes.}
\label{larger_final}
\begin{center}
\begin{footnotesize}
\begin{tabular}{c|c|r|rc|rc|rc|rc} \hline \hline
\rule{0mm}{1.2em} &   & & \multicolumn{2}{c|}{MILP} & \multicolumn{2}{c|}{MILP+MC} & \multicolumn{2}{c|}{CONIC} & \multicolumn{2}{c}{CONIC+MC} \\
\cline{3-11}
$\nu_0$ & $\kappa$ & \revtwo{assort} & \multicolumn{1}{c}{rgap} & \multicolumn{1}{c|}{time/\#} &\multicolumn{1}{c}{rgap} & \multicolumn{1}{c|}{time/\#} &\multicolumn{1}{c}{rgap} & \multicolumn{1}{c|}{time/\#} &\multicolumn{1}{c}{rgap} & \multicolumn{1}{c}{time/\#} \\
& & \revtwo{bind}& \multicolumn{1}{c}{egap} & \multicolumn{1}{c|}{nodes}& \multicolumn{1}{c}{egap} & \multicolumn{1}{c|}{nodes}& \multicolumn{1}{c}{egap} & \multicolumn{1}{c|}{nodes}& \multicolumn{1}{c}{egap} & \multicolumn{1}{c}{nodes} \\
\hline \hline
	&	20	& \revtwo{20.0} &	58.05\%	&	--	&	15.32\%	&	--	&	2.28\%	&	--	&	0.18\%	&	282.58/4	\\
	&		  & \revtwo{5}    &	57.54\%	&	188	&	14.73\%	&	114	&	0.37\%	&	1239	&	0.02\%	&	260	\\
\cline{2-11}	
	&	50	& \revtwo{50.0} &	32.14\%	&	--	&	11.14\%	&	--	&	5.56\%	&	--	&	0.11\%	&	188.05/5	\\
	&		  & \revtwo{5}    &	32.07\%	&	1235	&	11.05\%	&	546	&	2.61\%	&	1261	&	0.00\%	&	115	\\
	\cline{2-11}
10	&	100	& \revtwo{100.0} &	6.47\%	&	--	&	2.37\%	&	--	&	14.47\%	&	--	&	0.03\%	&	44.06/5	\\
	  &		  & \revtwo{5}     &	6.43\%	&	2022	&	2.18\%	&	3120	&	29.46\%	&	1371	&	0.00\%	&	6	\\
	\cline{2-11}
	  &	200	&	\revtwo{149.4} & 0.03\%	&	30.49/5	&	0.01\%	&	8.41/5	&	24.11\%	&	--	&	0.00\%	&	16.60/5	\\
	  &		  &	\revtwo{5} & 0.00\%	&	650	&	0.00\%	&	0	&	57.84\%	&	417	&	0.00\%	&	0	\\
	\cline{2-11}
	  &	500	&	\revtwo{149.4} & 0.03\%	&	38.30/5	&	0.02\%	&	13.04/5	&	24.11\%	&	--	&	0.00\%	&	18.30/5	\\
	  &		  & \revtwo{5}     &	0.00\%	&	756	&	0.00\%	&	10	&	55.57\%	&	64	&	0.00\%	&	0	\\
	\hline
	  &	20	&	\revtwo{20.0} & 24.48\%	&	--	&	9.57\%	&	--	&	1.35\%	&	--	&	0.04\%	&	165.95/5	\\
	  &		  &	\revtwo{5}    & 20.95\%	&	1109	&	9.41\%	&	127	&	0.10\%	&	1539	&	0.00\%	&	1	\\
	\cline{2-11}
		&	50	&	\revtwo{50.0} & 38.44\%	&	--	&	10.42\%	&	--	&	3.39\%	&	--	&	0.14\%	&	487.23/2	\\
		&			&	\revtwo{5}    & 38.44\%	&	840	&	10.37\%	&	322	&	0.75\%	&	1330	&	0.03\%	&	421	\\
\cline{2-11}
20	&	100	&	\revtwo{100.0} &	15.32\%	&	--	&	4.78\%	&	--	&	8.54\%	&	--	&	0.06\%	&	232.29/4	\\
		&			&	\revtwo{5}     & 15.30\%	&	1557	&	4.72\%	&	923	&	18.29\%	&	1430	&	0.01\%	&	276	\\
\cline{2-11}
	&	200	&	\revtwo{197.8} &	0.07\%	&	62.71/3	&	0.02\%	&	40.56/5	&	18.90\%	&	--	&	0.00\%	&	16.77/5	\\
	&			&	\revtwo{3}     &  0.06\%	&	7039	&	0.00\%	&	377	&	41.58\%	&	173	&	0.00\%	&	0	\\
\cline{2-11}
	&	500	&	\revtwo{203.4} & 0.02\%	&	15.84/5	&	0.01\%	&	9.31/5	&	19.90\%	&	--	&	0.00\%	&	18.37/5	\\
	&			&	\revtwo{0}     & 0.00\%	&	423	&	0.00\%	&	6	&	40.43\%	&	47	&	0.00\%	&	0	\\
\hline
\multicolumn{3}{c|}{Average}	&	17.51\%	&	33.96/18	&	5.36\%	&	17.83/20	&	12.26\%	&	--	&	0.06\%	&	119.43/45	\\
\multicolumn{3}{c|}{}			&	17.08\%	&	1582	&	5.25\%	&	555	&	24.70\%	&	887	&	0.01\%	&	108	\\
\hline \hline	
\end{tabular}
\end{footnotesize}
\end{center}
\end{table}
\egroup

\new{
 For the large instances, with the traditional (MILP) formulation the time limit is reached for 32 problem instances with tight capacity constraints.  Although the addition of McCormick inequalities reduces the integrality gaps substantially, only two more instances can be solved within time limit. The root gaps for the conic formulation are much smaller for the capacitated cases; nevertheless, problems cannot be solved to optimality within the time limit for these large instances. Adding the McCormick inequalities to the conic formulation reduces the average root gap to 0.06\% and allows the problems to be solved quickly. Many instances do not even require any branching, and 45 out of 50 instances are solved within the time limit. For the three instances that cannot be solved within the time limit, the end gap is only 0.04\% on average.
}

A third set of problems is inspired by the work of \citet{desir2014near}, who suggest a procedure to construct a family of hard benchmark instances \revtwo{to formally show that the MMNL assortment optimization problem is hard to approximate within any reasonable factor}. Each MMNL instance is generated based on an undirected graph $G=(V,E)$. Each vertex in $V$ corresponds to a product as well as a customer class ($V=M=N$). We denote by $C_i=\{ j \,| \,(i,j) \in E \}$ the set of products that the customers in class $i$ consider buying (this always includes product $i$ and can be thought of as class $i$'s ``consideration set'').  Given this structure, we create a problem set with 100 products (and 100 classes). Each product has 10 neighbors in $G$ so  $|C_i|=11$. These neighbors are selected randomly. \revtwo{However, this} procedure may lead to unrealistic preference and price parameters; \revtwo{therefore,} we use the following \revtwo{modification}. We denote product $i$ as class $i$'s favorite product and set $\nu_{ii}=1$.  For $i \neq j, (i,j) \in E$, $\nu_{ij}$ is drawn from a uniform$[0,1]$ distribution. For  $(i,j) \notin E$, $\nu_{ij}=0$. The prices are randomly generated from a uniform$[1,3]$ distribution.  $\gamma_i, i \in M$ are drawn from a uniform$[0,1]$ distribution. The capacity $\kappa$ is one of $\{10,20,50,100\}$, the no-purchase parameter $\nu_{i0}$ is either 1 or 2, and we again generate five instances for each parameter setting, leading to 40 instances. The results are reported in Table \ref{harder_final}.

\bgroup
\setlength\tabcolsep{0.08cm}
\def\arraystretch{1}%
\begin{table}[h]
\caption{Results for hard problems.}
\label{harder_final}
\begin{center}
\begin{footnotesize}
\begin{tabular}{c|c|r|rc|rc|rc|rc} \hline \hline
\rule{0mm}{1.2em} &  &  & \multicolumn{2}{c|}{MILP} & \multicolumn{2}{c|}{MILP+MC} & \multicolumn{2}{c|}{CONIC} & \multicolumn{2}{c}{CONIC+MC} \\
\cline{3-11}
$\nu_0$ & $\kappa$ & \revtwo{assort}  &  \multicolumn{1}{c}{rgap} & \multicolumn{1}{c|}{time/\#} &\multicolumn{1}{c}{rgap} & \multicolumn{1}{c|}{time/\#} &\multicolumn{1}{c}{rgap} & \multicolumn{1}{c|}{time/\#} &\multicolumn{1}{c}{rgap} & \multicolumn{1}{c}{time/\#} \\
& & \revtwo{bind} & \multicolumn{1}{c}{egap} & \multicolumn{1}{c|}{nodes}& \multicolumn{1}{c}{egap} & \multicolumn{1}{c|}{nodes}& \multicolumn{1}{c}{egap} & \multicolumn{1}{c|}{nodes}& \multicolumn{1}{c}{egap} & \multicolumn{1}{c}{nodes} \\
\hline \hline
 	&	10	& \revtwo{10.0} &	32.60\%	&	17.29/5	&	6.50\%	&	9.30/5	&	8.28\%	&	209.88/1	&	1.75\%	&	4.22/5	\\
	&		  & \revtwo{5} &	0.00\%	&	2662	&	0.00\%	&	1154	&	1.04\%	&	15614	&	0.00\%	&	84	\\
\cline{2-11}	
 	&	20	& \revtwo{20.0} &	27.93\%	&	--	&	10.99\%	&	--	&	8.19\%	&	349.10/2	&	1.71\%	&	14.20/5	\\
1	&		 	& \revtwo{5} &	4.07\%	&	20345	&	2.04\%	&	28002	&	0.98\%	&	8680	&	0.00\%	&	370	\\
\cline{2-11}	
 	&	50	& \revtwo{50.0} &	3.15\%	&	89.21/5	&	0.80\%	&	3.47/5	&	12.96\%	&	--	&	0.12\%	&	1.22/5	\\
	&		  & \revtwo{5}    &	0.00\%	&	14380	&	0.00\%	&	1032	&	2.02\%	&	9343	&	0.00\%	&	0	\\
\cline{2-11}	
	&	100	& \revtwo{64.2} &	1.25\%	&	3.94/5	&	0.21\%	&	0.49/5	&	13.57\%	&	--	&	0.06\%	&	0.50/5	\\
	&		  & \revtwo{0}    &	0.00\%	&	6099	&	0.00\%	&	23	&	1.20\%	&	10950	&	0.00\%	&	0	\\
\hline
	&	10	& \revtwo{10.0} &	12.90\%	&	7.93/5	&	2.92\%	&	4.24/5	&	3.79\%	&	124.74/5	&	0.61\%	&	2.61/5	\\
	&		  & \revtwo{5}    &	0.00\%	&	1565	&	0.00\%	&	484	&	0.00\%	&	7514	&	0.00\%	&	38	\\
\cline{2-11}	
 	&	20	& \revtwo{20.0} &	21.82\%	&	298.66/1	&	6.26\%	&	267.82/2	&	4.55\%	&	129.21/4	&	0.69\%	&	8.62/5	\\
2	&		  & \revtwo{5}    &	0.75\%	&	42703	&	0.65\%	&	32411	&	0.04\%	&	6233	&	0.00\%	&	188	\\
\cline{2-11}	
 	&	50	& \revtwo{50.0} &	6.24\%	&	482.84/3	&	1.30\%	&	47.18/4	&	8.06\%	&	--	&	0.20\%	&	6.64/5	\\
	&		  & \revtwo{5} &	0.26\%	&	23303	&	0.08\%	&	20401	&	0.81\%	&	10649	&	0.00\%	&	251	\\
\cline{2-11}	
 	&	100	& \revtwo{79.8} &	0.39\%	&	1.01/5	&	0.01\%	&	0.18/5	&	8.01\%	&	183.98/4	&	0.00\%	&	0.39/5	\\
	&		  & \revtwo{0}    &	0.00\%	&	445	&	0.00\%	&	0	&	0.07\%	&	6950	&	0.00\%	&	0	\\
\hline
 \multicolumn{3}{c|}{Average}		 	&	13.29\%	&	80.83/29	&	3.63\%	&	26.22/31	&	8.43\%	&	174.03/16	&	0.64\%	&	4.8/40	\\
 \multicolumn{3}{c|}{}		&	0.64\%	&	13938	&	0.35\%	&	10438	&	0.77\%	&	9492	&	0.00\%	&	116	\\
\hline \hline	
\end{tabular}
\end{footnotesize}
\end{center}
\end{table}
\egroup

\new{
These instances are indeed harder than the previous sets.
The root gaps for the (CONIC+MC) formulation are higher than those of the previous sets. Nevertheless,
the relative effectiveness of the formulations is consistent with the earlier experiments.
With the (CONIC+MC) formulation, all instances are solved within the time limit with an average run time under five seconds.
}

\new{
In the final set of experiments, we compare the formulations on instances with generalized capacity constraints. The general capacity data set has 200 products and 20 classes. The preference values and prices are generated as in Table \ref{smaller_final}. The model has six capacity constraints. The first constraint is a general capacity constraint $\sum_{j \in S} \beta_{0j} \leq \kappa_0$,
where $\beta_{0j}$ is generated uniformly between 0 and 1.
The other five constraints are subset cardinality constraints $|S \cap N_{k}| \leq \kappa_k$, $k=1,\ldots,5$ where $N_{k}$, $k=1,\ldots,5$ are disjoint sets with $|N_{k}|=40$. In order to obtain the lower bounds for the conditional McCormick inequalities, we use the following approach:
For both conditions ($x_j=1$ and $x_j=0$), we first solve the linear relaxation of (BND) with only the capacity constraint using the greedy algorithm. We then solve the same problem with only the non-overlapping subset cardinality constraints \revtwo{also using the greedy algorithm}. We use the minimum of the two relaxation values to obtain the lower bounds. \revtwo{Separately considering the constraints allows us to utilize fast greedy algorithms instead of using the simplex or an interior point algorithm for each variable-value combination.} The results are shown in Table \ref{generalized_capacity}, \revtwo{where space reports the amount of capacity used ($\sum_{j \in N} \beta_{0j}x_j^*$, averaged over five instances) and bind now reports the number of instances where {\it all} subset cardinality constraints are tight in the optimal solution}.
}

\bgroup
\setlength\tabcolsep{0.08cm}
\def\arraystretch{1}%
\begin{table}[h]
\caption{Results for problems with generalized capacity constraints.}
\label{generalized_capacity}
\begin{center}
\begin{footnotesize}
\begin{tabular}{c|c|r|rc|rc|rc|rc} \hline \hline
\rule{0mm}{1.2em} &  & & \multicolumn{2}{c|}{MILP} & \multicolumn{2}{c|}{MILP+MC} & \multicolumn{2}{c|}{CONIC} & \multicolumn{2}{c}{CONIC+MC} \\
\cline{4-11}
$\nu_0$ & $\kappa_0, \kappa_k$ & \revtwo{space/assort}&  \multicolumn{1}{c}{rgap} & \multicolumn{1}{c|}{time/\#} &\multicolumn{1}{c}{rgap} & \multicolumn{1}{c|}{time/\#} &\multicolumn{1}{c}{rgap} & \multicolumn{1}{c|}{time/\#} &\multicolumn{1}{c}{rgap} & \multicolumn{1}{c}{time/\#} \\
& & \revtwo{bind} & \multicolumn{1}{c}{egap} & \multicolumn{1}{c|}{nodes}& \multicolumn{1}{c}{egap} & \multicolumn{1}{c|}{nodes}& \multicolumn{1}{c}{egap} & \multicolumn{1}{c|}{nodes}& \multicolumn{1}{c}{egap} & \multicolumn{1}{c}{nodes} \\
\hline \hline
 	&	5,2 & \revtwo{4.48/10.0}	&	24.80\%	&	--	&	7.65\%	&	--	&	1.95\%	&	194.98/4	&	0.10\%	&	5.97/5	\\
	&		  & \revtwo{5}          &	10.07\%	&	53576	&	5.51\%	&	6836	&	0.11\%	&	9351	&	0.00\%	&	4	\\
\cline{2-11}	
 	&	10,4 & \revtwo{9.42/20.0} &	40.82\%	&	--	&	9.32\%	&	--	&	3.78\%	&	--	&	0.25\%	&	12.12/5	\\
	&		   &          \revtwo{5} &	38.73\%	&	3414	&	7.99\%	&	10972	&	1.48\%	&	7870	&	0.00\%	&	97	\\
\cline{2-11}
10	& 25,10 & \revtwo{24.11/50.0} &	13.16\%	&	--	&	3.84\%	&	--	&	10.02\%	&	--	&	0.44\%	&	114.38/5	\\
	  &       &           \revtwo{5} &	11.89\%	&	49065	&	2.68\%	&	36783	&	4.22\%	&	3008	&	0.00\%	&	1971	\\
\cline{2-11}
	&	50,20	&	\revtwo{45.19/87.6} & 0.12\%	&	4.81/5	&	0.04\%	&	1.47/5	&	17.00\%	&	--	&	0.01\%	&	2.56/5	\\
	&		    &	           \revtwo{0} & 0.00\%	&	1291	&	0.00\%	&	27	&	5.50\%	&	767	&	0.00\%	&	26	\\
\cline{2-11}
 	&	100,40	& \revtwo{45.66/88.6} &	0.04\%	&	1.17/5	&	0.02\%	&	0.9/5	&	17.19\%	&	--	&	0.00\%	&	1.81/5	\\
	&		      &	\revtwo{0} & 0.00\%	&	95	&	0.00\%	&	0	&	4.47\%	&	1406	&	0.00\%	&	0	\\
\hline
	&	5,2	  &	\revtwo{4.50/10.0} & 6.94\%	&	--	&	2.61\%	&	479.15/1	&	1.03\%	&	48.86/5	&	0.02\%	&	4.72/5	\\
	&		    &	\revtwo{5}          & 1.01\%	&	136413	&	0.63\%	&	32437	&	0.00\%	&	3120	&	0.00\%	&	0	\\
\cline{2-11}
 	&	10,4	  &	\revtwo{9.36/20.0} & 21.43\%	&	--	&	5.25\%	&	--	&	2.15\%	&	--	&	0.10\%	&	10.75/5	\\
	&		      & \revtwo{5}          &	14.77\%	&	21448	&	4.35\%	&	9093	&	0.63\%	&	15089	&	0.00\%	&	55	\\
\cline{2-11}
20 &	25,10	&	\revtwo{24.28/50.0} & 18.45\%	&	--	&	3.87\%	&	--	&	5.73\%	&	--	&	0.24\%	&	97.01/5	\\
	 &		    & \revtwo{5}           & 17.73\%	&	12587	&	3.10\%	&	21435	&	2.07\%	&	4784	&	0.00\%	&	1559	\\
\cline{2-11}
	&	50,20	&	\revtwo{49.20/98.6} & 1.54\%	&	206.11/2	&	0.32\%	&	21.57/4	&	11.66\%	&	--	&	0.06\%	&	10.58/5	\\
	&		    & \revtwo{2}           &	0.79\%	&	135867	&	0.04\%	&	9802	&	3.37\%	&	868	&	0.00\%	&	122	\\
\cline{2-11}
	&	100,40 & \revtwo{60.81/120.2} & 0.03\% &	0.99/5	&	0.01\%	&	0.75/5	&	12.82\%	&	--	&	0.00\%	&	1.41/5	\\
	&				 & \revtwo{0}           &	0.00\% &	99	&	0.00\%	&	0	&	5.20\%	&	1321	&	0.00\%	&	0	\\
\hline	
 \multicolumn{3}{c|}{Average}		&	12.73\%	&	26.30/17	&	3.29\%	&	29.05/20	&	8.33\%	&	113.80/9	&	0.12\%	&	26.13/50	\\
 \multicolumn{3}{c|}{}		&	9.50\%	&	41385	&	2.43\%	&	12738	&	2.71\%	&	4758	&	0.00\%	&	383	\\
\hline \hline	
\end{tabular}
\end{footnotesize}
\end{center}
\end{table}
\egroup

\new{
The results in Table \ref{generalized_capacity} are consistent with earlier experiments. The CONIC+MC formulation leads to tight relaxations under generalized capacity constraints as well. All 50 instances are solved in under 30 seconds on average, whereas with the second best formulation (MILP+MC), only 20 instances are solved. We note that the time to compute the conditional bounds is negligible as we utilize a greedy approach to solve the relaxations.
}

\section{Concluding Remarks}
In this paper, we present a conic quadratic mixed-integer formulation of the capacitated assortment optimization problem under the mixed multinomial logit model, that is far more effective than traditional MILP formulations of this problem with tight capacity constraints. Additional performance improvements are gained by using McCormick estimators derived through conditional bounds exploiting the capacity constraints. The numerical results suggest that with the new formulations, commercially available software may be practically used to solve even relatively large assortment optimization problems to optimality. Given the promise of conic mixed-integer formulations for the MMNL problem, it is worthwhile to explore conic optimization formulations of assortment optimization problems based on other consumer choice models.

\section*{Acknowledgements}
A. \c{S}en was supported by a 2219 fellowship grant from The Scientific and Technological Research Council of Turkey (T\"{U}B\.{I}TAK).
He acknowledges with gratitude the financial support of T\"{U}B\.{I}TAK and hospitality of University of California-Berkeley.
A. Atamt\"urk was supported, in part, by grant FA9550-10-1-0168 from the Office of
the Assistant Secretary of Defense for Research and Engineering.
P. Kaminsky was supported, in part, by industry members of the I/UCRC
Center for Excellence in Logistics and Distribution, and by
National Science Foundation Grant No. 1067994.


\bibliographystyle{apalike}
\bibliography{AssortmentBib}

\end{document}